\documentclass[12pt,oneside]{amsart}
\usepackage{geometry}                
\usepackage{graphicx}
\usepackage{amssymb}
\usepackage{epstopdf}
\usepackage{amsmath,amsthm,amscd,amssymb}
\usepackage{latexsym}
\usepackage[colorlinks,citecolor=red,pagebackref,hypertexnames=false]{hyperref}

\numberwithin{equation}{section}

\theoremstyle{plain}
\newtheorem{theorem}{Theorem}[section]
\newtheorem{lemma}[theorem]{Lemma}
\newtheorem{corollary}[theorem]{Corollary}
\newtheorem{proposition}[theorem]{Proposition}

\theoremstyle{definition}

\theoremstyle{remark}

\newtheorem{case[theorem]}{Case}

\date{March 4, 2011}      
\author{Alex Iosevich and Krystal Taylor}
\address{Department of Mathematics, University of Rochester, Rochester, NY}
\email{iosevich@math.rochester.edu}
\address{Department of Mathematocs, University of Rochester, Rochester, NY}
\email{taylor@math.rochester.edu}

\thanks{This work was partially supported by the NSF Grant DMS10-45404.}

\title{\parbox{14cm}{\centering{Lattice points close to families of surfaces, non-isotropic dilations and regularity of generalized Radon transforms}}}

\begin{document}
\maketitle


\begin{abstract} We prove that if $\phi: {\Bbb R}^d \times {\Bbb R}^d \to {\Bbb R}$, $d \ge 2$, is a homogeneous function, smooth away from the origin and having  non-zero Monge-Ampere determinant away from the origin, then 
$$ R^{-d} \# \{(n,m) \in {\Bbb Z}^d \times {\Bbb Z}^d: |n|, |m| \leq CR; R \leq \phi(n,m) \leq R+\delta \} \lesssim \max \{R^{d-2+\frac{2}{d+1}}, R^{d-1} \delta \}.$$ 

This is a variable coefficient version of a result proved by Lettington in \cite{L10}, extending a previous result by Andrews in \cite{A63}, showing that if $B \subset {\Bbb R}^d$, $d \ge 2$, is a symmetric convex body with a sufficiently smooth boundary and non-vanishing Gaussian curvature,  then 
$$ \# \{k \in {\mathbb Z}^d: dist(k, R \partial B) \leq \delta \} \lesssim \max \{R^{d-2+\frac{2}{d+1}}, R^{d-1} \delta \}. (*)$$ 

Furthermore, we shall see that the same argument yields a non-isotropic analog of $(*)$, one for which the exponent on the right hand side is, in general, sharp, even in the infinitely smooth case. This sheds some light on the nature of the exponents and their connection with the conjecture due to Wolfgang Schmidt on the distribution of lattice points on dilates of smooth convex surfaces in ${\Bbb R}^d$. 
\end{abstract}  

\maketitle


\section{Introduction} 

The problem of counting integer lattice points inside, on, and near convex surfaces is a classical and time-honored problem in number theory and related areas. See \cite{Hux96} and the references contained therein for a thorough description of this beautiful area. In this paper we shall focus on the problem of counting integer lattice points in the neighborhood of variable coefficient families of surfaces. It follows from a result of G. Andrews (\cite{A63}) that if $B \subset {\Bbb R}^d$, $d \ge 2$, is a symmetric convex body with a strictly convex boundary, then 
\begin{equation} \label{andrews} \# \{R \partial B \cap {\Bbb Z}^d \} \lesssim R^{d-2+\frac{2}{d+1}}, \end{equation} where the implicit constant depends on $B$ and the dimension. 

It is not known to what extent (\ref{andrews}) is sharp, at least in higher dimensions. In dimension two, one can show that there exists an infinite sequence of $R$s going to infinity, such that 
\begin{equation} \label{2bad}  \# \{R \partial B \cap {\Bbb Z}^2 \} \gtrsim C_{\epsilon} R^{\frac{2}{3}-\epsilon}\end{equation} for any $\epsilon>0$. See, for example, \cite{P06} and \cite{ISS07}. It is important to note, however, that the boundary of $B$ in (\ref{2bad}) is only $C^{1,1}$ and not any smoother. 

In dimensions three and higher, a deep and far-reaching conjecture due to Wolfgang Schmidt (\cite{Sch86}) says that if the boundary of $B$ is smooth and has non-vanishing Gaussian curvature, then for any 
$\epsilon>0$,
\begin{equation} \label{schmidt} \# \{R \partial B \cap {\Bbb Z}^d \} \leq C_{\epsilon} R^{d-2+\epsilon}. \end{equation} 

See \cite{IR07} for the discussion of related issues. In dimension two, even smoothness does not lead to an analog of (\ref{schmidt}), even conjecturally, due to an example due to Konyagin (\cite{K77}), who showed that there exists a smooth symmetric convex curve $\Gamma$, with everywhere non-vanishing curvature, and a sequence of dilates $R_j \to \infty$ such that 
$$\# \{R_j \Gamma \cap {\Bbb Z}^2 \} \gtrsim \sqrt{R_j}.$$ 

\vskip.125in 

M. C. Lettington (\cite{L10}) recently extended Andrews' result ((\ref{andrews}) above) by showing that 
\begin{equation} \label{key} \# \{k \in {\Bbb Z}^d: R \leq {||k||}_B \leq R+\delta \} \lesssim \max \{ R^{d-2+\frac{2}{d+1}}, R^{d-1} \delta \}, \end{equation}
where 
$${||x||}_B=\inf \{t>0: x \in tB \}.$$ He needs $\partial B$ to have a tangent hyper-plane at every point and that any two-dimensional cross section through the normal consist of a plane curve with continuous radius of curvature bounded away from zero and infinity. 

If the boundary of $B$ is smoother, Lettington's bound can be improved in the following way. Let 
$$ N_B(R)=\# \{RB \cap {\Bbb Z}^d \}. $$

Define the discrepancy function, ${\mathcal D}_B(R)$ by the equation 
$$ N_B(R)=|B|R^d+{\mathcal D}_B(R).$$ 

Suppose that we have a bound 
\begin{equation} \label{discrepancy} |{\mathcal D}_B(R)| \lessapprox R^{d-2+\alpha_d} \end{equation} for some $\alpha_d>0$. \footnote{Here and throughout, $X \lessapprox Y$, with the controlling parameter $R$, if for every $\epsilon>0$ there exists $C_{\epsilon}>0$ such that $X \leq C_{\epsilon} R^{\epsilon} Y$}. It follows that 
$$ |N_B(R+\delta)-N_B(R)|= \left| |B|{(R+\delta)}^d+{\mathcal D}_B(R+\delta)-|B|R^d-{\mathcal D}_B(R) \right|$$
$$ \leq \left||B|{(R+\delta)}^d-|B|R^d \right|+|{\mathcal D}_B(R+\delta)|+|{\mathcal D}_B(R)|$$
\begin{equation} \label{discrepancyreacharound} \lessapprox R^{d-1} \delta+R^{d-2+\alpha_d},\end{equation}  which is an immediate improvement over (\ref{key}) if $\alpha_d<\frac{2}{d+1}$. 

Indeed, Wolfgang Muller (\cite{Mu99}) proves that (\ref{discrepancy}) holds with 
$$ \alpha_d=\frac{d+4}{d^2+d+2} \ \text{if} \ d \ge 5; \ \alpha_4=\frac{6}{17} \ \text{and} \ \alpha_3=\frac{20}{43}.$$ 

It is not difficult to check that in each case, $\alpha_d<\frac{2}{d+1}$. 

The purpose of this paper is two-fold. We extend Lettington's estimate to a variable coefficient setting where generalized Radon transforms play the dominant role. In the process, we give a reasonably short Fourier analytic proof of (\ref{key}) under more stringent smoothness assumptions on $\partial B$ than the ones used by Lettington, but less stringent than those needed by Muller.  We shall also see that the same argument yields a certain multi-parameter analog of (\ref{key}), one for which the exponent $d-2+\frac{2}{d+1}$ in (\ref{key}) \underline{cannot} be improved, even in the {\it infinitely smooth case}. We shall also obtain a non-isotropic variant of (\ref{discrepancy}) where, once again, the exponent $d-2+\frac{2}{d+1}$ cannot be improved, even in the infinitely smooth case. This sheds some light on the nature of the exponents and further illustrates the depth of Schmidt's conjecture (\ref{schmidt}). Our main results, initially stated in an isotropic setting, are the following. 

\begin{theorem} \label{main} Let $\phi: {\Bbb R}^d \times {\Bbb R}^d \to {\Bbb R}$ be a homogeneous function of degree one, $C^{\lfloor \frac{d}{2} \rfloor+1}$ away from the origin. Suppose that 
\begin{equation} \label{manifoldstructure} \nabla_x \phi_l(x,y)\neq \vec{0} \text{ and } \nabla_y \phi_l(x,y)\neq \vec{0} \end{equation} in a neighborhood of the sets 
$$\{x \in B: \phi(x,y)=t\}, \{y \in B: \phi(x,y)=t\},$$ where, $B$ denotes the unit ball. Suppose further that the Monge-Ampere determinant of $\phi$, (introduced by Phong and Stein the Euclidean setting in \cite{PS86}), given by 
\begin{equation} \label{mongeampere} \det \begin{pmatrix} 
 0 & \nabla_{x}\phi \\
 -{(\nabla_{y}\phi)}^{T} & \frac{\partial^2 \phi}{dx_i dy_j}
\end{pmatrix}, \end{equation} does not vanish on the set $\{(x,y) \in B \times B: \phi(x,y)=t \}$ for any $t>0$. 

Then 

\begin{equation} \label{maineq} q^{-d} \# \{(n,m) \in {\Bbb Z}^d \times {\Bbb Z}^d: |n|, |m| \leq Cq; q \leq \phi(n,m) \leq q+\delta \} \lesssim max \{q^{d-2+\frac{2}{d+1}}, q^{d-1} \delta \}, \end{equation}

for $C$ a positive constant dependent on $\phi$.   
\end{theorem} 

\begin{corollary} \label{mainbonus} Let $B$ be a bounded symmetric convex body. Suppose that $\partial B$ is $C^{\lfloor \frac{d}{2} \rfloor+1}$ and has everywhere non-vanishing Gaussian curvature. Then (\ref{key}) holds. \end{corollary} 

Corollary \ref{mainbonus} follows from Theorem \ref{main} by first observing that if $\partial B$ is $C^{\lfloor\frac{d}{2} \rfloor+1}$ and has everywhere non-vanishing Gaussian curvature, then $\phi(x,y)={||x-y||}_B$ satisfies the Monge-Ampere assumption in (\ref{mongeampere}) above, as can be demonstrated by a direct calculation. This gives us (\ref{maineq}). We then observe that if $\phi(x,y)={||x-y||}_B$ and $R=q$, the left hand side in (\ref{maineq}) equals the left hand side of (\ref{key}). This completes the proof of the corollary, assuming Theorem \ref{main}.

\subsection{Non-isotropic formulation and sharpness of exponents} 

We are going to prove the following, more general, version of Theorem \ref{main}. 

\begin{theorem} \label{maingeneral} Let $\phi: {\Bbb R}^d \times {\Bbb R}^d \to {\Bbb R}$ be a
$C^{\lfloor \frac{d}{2} \rfloor+1}$ function away from the origin satisfying the quasi-homogeneity condition
$$\phi(q^{\alpha_1}x_1,\cdots, q^{\alpha_d}x_d, q^{\alpha_1}y_1,..., q^{\alpha_d}y_d)=q^{\beta}\phi(x,y),$$ where $\sum_{j=1}^d\alpha_j=d$, $\alpha_j \le \frac{2d}{d+1}$, and $\alpha_j, \beta>0$.
Suppose further that the Monge-Ampere determinant of $\phi$, given in (\ref{mongeampere}), does not vanish for any $t>0$, and that the non-degeneracy assumption (\ref{manifoldstructure}) holds. Then 

\begin{equation} \label{generaleq}  q^{-d}\# \{(n,m) \in {\Bbb Z}^d \times {\Bbb Z}^d: \forall j, |n_j|, |m_j| \leq Cq^{\alpha_j}; | \phi(n,m)- q^{\beta}|\leq \delta \}
\end{equation} 
$$\lesssim max \{q^{d-2+\frac{2}{d+1}}, q^{d-\beta} \delta \}, $$
for a positive constant $C$ dependent on $\phi$.
\end{theorem}

Theorem \ref{main} follows from Theorem \ref{maingeneral} by taking $\beta=\alpha_j \equiv 1$. 

Once again, in the case when $\phi(n,m)=\phi_0(n-m)$, we get a non-isotropic analog of Corollary \ref{mainbonus}. 

\subsubsection{Sharpness of exponents:} To see that Theorem \ref{maingeneral} is, in general, sharp, let 
$$\phi(x,y)=(x_d-y_d)-{(x_1-y_1)}^2-\dots-{(x_{d-1}-y_{d-1})}^2,$$ 
\begin{equation} \label{coolalpha} \alpha_1=\dots=\alpha_{d-1}=\frac{d}{d+1}, \ \alpha_d=\beta=\frac{2d}{d+1}. \end{equation} 

It is not hard to see that with this $\phi$,

$$ q^{-d} \# \{(n,m) \in {\Bbb Z}^d \times {\Bbb Z}^d:\forall j, |n_j|, |m_j| \leq Cq^{\alpha_j}; \phi(n,m)=q^{\frac{2d}{d+1}} \} \approx q^{d-2+\frac{2}{d+1}},$$ and thus Theorem \ref{maingeneral} is sharp. We note that in a discrete two-dimensional setting, this type of a construction was used by Pavel Valtr (\cite{V05}) to give an example of a family of points and translates of a fixed convex curve with everywhere non-vanishing curvature for which the exponent given by the Szemeredi-Trotter incidence theorem cannot be improved. See also \cite{IS10} where Valtr example is explored in a continuous setting of the Falconer distance conjecture. 

Going back to Wolfgang Schmidt's conjecture (\ref{schmidt}), we see that our example above clearly shows that isotropic dilations are absolutely necessary for the conjecture to hold. Using non-isotropic dilations, the conjectured exponent $d-2$ may be as bad as $d-2+\frac{2}{d+1}$. These observations suggests that there is a delicate interplay between the smoothness of the boundary and the structure of dilations which should prove to be a fruitful field of investigation in the sequel. 

\vskip.125in 

{\bf Acknowledgements:} The authors are grateful to D. H. Phong for some helpful remarks on this paper.

\section{Proof of the main result (Theorem (\ref{maingeneral}))}

The argument below is motivated, to a significant degree, by Falconer's argument in \cite{Fal86}. See also \cite{Falc86} and \cite{M95}. See \cite{IJL10} and \cite{EIT10} where Sobolev bounds for generalized Radon transforms are used to obtain geometric and geometric combinatorial conclusions. 

Set $$\mu_q(x)= q^{-d} q^{\frac{d^2}{s}} \sum_{a \in \mathbb{Z}^d} \prod_{j=1}^d \psi_0 \left(\frac{a_j}{q^{\alpha_j}}\right) \psi_0 \left(q^{\frac{d}{s}}\left(x_j-\frac{a_j}{q^{\alpha_j}}\right)\right),$$ where $\psi_0$ is a smooth symmetric function which is identically equal to $1$ on the unit ball and equal to $0$ outside of the interval $(-C,C)$ for a positive constant $C>2$. 

Let $E_q$ denote the support of $\mu_q$. Notice that $\mu_q(B(x, q^{-\frac{d}{s}})) \sim q^{-d}$ for $x\in E_q$.

We show that
\begin{equation} \label{keyestimate} 
\mu_q \times \mu_q (\{ (x,y): |\phi(x,y)-1| \le q^{-\frac{d}{s}} \}) \lesssim q^{-\frac{d}{s}} \end{equation}
for $\frac{d+1}{2}\le s < d$. 
\\

It is an immediate consequence of this estimate that
$$q^{-d}\# \{(n,m) \in {\Bbb Z}^d \times {\Bbb Z}^d: \forall j, |n_j|, |m_j| \leq Cq^{\alpha_j}; |\phi(n,m)- q^{\beta}| \le \delta \} \lesssim max \{q^{d-2+\frac{2}{d+1}}, q^{d-\beta} \delta \}. $$ 

Indeed, letting $N(A, \gamma)$ denotes the number of balls of radius $\gamma$ needed to cover a set $A$ and letting $\tau^{\alpha}_qx=(q^{\alpha_1}x_1, \dots, q^{\alpha_d} x_d)$, we get
\begin{align*}
&\mu_q \times \mu_q (\{ (x,y): |\phi(x,y)-1| \le q^{-\frac{d}{s}}\}\\
&\sim q^{-2d} N(\{ (x,y)\in E_q \times E_q: |\phi(x,y)-1| \le q^{-\frac{d}{s}}\}, q^{-\frac{d}{s}} )\\
&\sim q^{-2d} N(\{ (u,v)\in \tau^{\alpha}_q(E_q) \times \tau^{\alpha}_q(E_q): |\phi(u,v)-q^{\beta}| \le q^{\beta-\frac{d}{s}}\}, q^{\beta-\frac{d}{s}})\\
&\sim q^{-2d} \#(\{ (n,m)\in {\Bbb Z}^d \times {\Bbb Z}^d: \forall j, |n_j|, |m_j| \leq Cq^{\alpha_j}, |\phi(n,m)-q^{\beta}| \leq q^{\beta-\frac{d}{s}}\}).\\
\end{align*} 

If $\delta \le q^{\beta- \frac{2d}{d+1}}$ then result follows immediatley by setting $s=\frac{d+1}{2}$. 
Otherwise, choose $\frac{d+1}{2}\le s <d$ so that $\delta = q^{\beta-\frac{d}{s}}$.   
\\

To show (\ref{keyestimate}), we begin by re-writing the left hand side of the inequality as
$$\mu_q \times \mu_q (\{ (x,y): |\phi(x,y)-1| \le q^{-\frac{d}{s}}\})
=\iint_{\{|\phi(x,y)-1| \le q^{-\frac{d}{s}} \}}\psi(x,y) d\mu_q(y)d\mu_q(x)$$
where $\psi$ is a smooth function with compact support which is centered at the origin.

Define
\begin{equation}\label{operator} T_q f(x)= q^{\frac{d}{s}}\int_{\{|\phi(x,y)-1| \le q^{-\frac{d}{s}} \}}f(y)\psi(x,y) dy.\end{equation}

Then the left hand side of $(\ref{keyestimate})$ can be written as $ <T_q\mu_q ,\mu_q >$, and it remains to show that, for $\frac{d+1}{2}\le s < d$,
$$ <T_q\mu_q ,\mu_q >\lesssim 1.$$ 

By  the Cauchy Schwarz inequality:
\begin{equation} \label{heaven}  <T_q\mu_q ,\mu_q > \le \|(\widehat{T_q\mu_q}(\cdot)\times|\cdot|^{\frac{d-s}{2}}) \|_2  \times \|(\widehat{\mu_q}(\cdot)|\cdot|^{\frac{s-d}{2}}) \|_2=I \times II. \end{equation} 

The remainder of the paper is dedicated to showing that terms $I$ and $II$ are bounded. As we point out above, this implies (\ref{keyestimate}) which, in turn, implies Theorem \ref{maingeneral}. 

\subsection{Estimation of the first term in (\ref{heaven})} 

\begin{lemma}\label{energy} 
Denote the $s-dimensional$ energy of $\mu_q$ by $I_s(\mu_q)$ (see, for instance, \cite{W04}). Then $$I_s(\mu_q)=\|(\widehat{\mu_q}(\cdot)|\cdot|^{\frac{s-d}{2}}) \|_2^2 \lesssim 1.$$ 
\end{lemma}

To prove the lemma, observe that 
\begin{equation}\label{acro}
\int|\widehat{\mu}_q(\xi)|^2 |\xi|^{s-d}  d\xi 
=\iint|x-y|^{-s} d\mu_q(x) d\mu_q(y).
\end{equation}

We expand the right hand side of (\ref{acro}) using the definition of $\mu_q$. In doing so, we introduce the summation over $a \in \mathbb{Z}^d$ and $a' \in \mathbb{Z}^d$. 

\vskip.125in

{\bf The isotropic case:} To motivate the argument for the general (non-isotropic) case which follows, we first look at the proof of this lemma in the isotropic case, the regime where $\alpha_j=1$ for all $1 \le j \le d$. In this case, (\ref{acro}) becomes

$$q^{-2d}q^{\frac{2d^2}{s}}
\sum_{a, a' \in \mathbb{Z}^d}  \psi_0 \left(\frac{a}{q}\right)  \psi_0 \left(\frac{a'}{q}\right) 
\iint \psi_0 \left(q^{\frac{d}{s}}\left(x -\frac{a}{q}\right)\right)\psi_0 \left(q^{\frac{d}{s}}\left(y-\frac{a'}{q}\right)\right)|x-y|^{-s}dxdy. $$
When $a=a'$,  the above quantity is certainly bounded. Indeed, if $a=a'$ then both $x$ and $y$ lie the same ball of radius $\sim q^{-\frac{d}{s}}$, and integrating with spherical coordinates gives the desired result. 

In the case that $a\neq a'$, we reduce our problem of bounding (\ref{acro}) to bounding 
\begin{equation}\label{simplecase} q^{s}q^{-2d} \sum_{a, a' \in \mathbb{Z}^d}  \psi_0 \left(\frac{a}{q}\right)  \psi_0 \left(\frac{a'}{q}\right) |a-a'|^{-s}   .\end{equation}

To accomplish this, we break the sum into dyadic shells. For a fixed  $a' \in \mathbb{Z}^d$ with $|a'|\le C q$, set 
$$A_m = \{ a \in \mathbb{Z}^d: 2^m \le |a-a'| < 2^{m+1} \}$$ where $0\le m \le \log{ \lceil Cq\rceil}$. 
We use the fact that $\#(A_m) \sim 2^{md}$. 

Now 
$$ \sum_{a \in \mathbb{Z}^d}  \psi_0 \left(\frac{a}{q}\right)  |a-a'|^{-s}   
\sim \sum_m \sum_{a \in A_m}    |a-a'|^{-s}   
\lesssim  \sum_m 2^{m(d-s)}
\sim q^{(d-s)}.$$
Plugging this calculation into (\ref{simplecase}) and summing in   $a' \in \mathbb{Z}^d$ with $|a'|\le C q$, we see that $(\ref{acro})\lesssim 1$. 

\vskip.125in 

\textbf{The general (non-isotropic) case in two dimensions:} For greater clarity, before we proceed to the non-isotropic case for general $d$, we look at the specific case when $d=2$.  We again break the sum over $a, a' \in \mathbb{Z}^2$ into the sum over the diagonal and the sum away from the diagonal. 
When  $a=a'$, we use the same technique as above. 

In the case that $a\neq a'$, set 
$$I(a,a') = \iint  \prod_{j=1}^2 \psi_0 \left(q^{\frac{2}{s}}\left(x_j-\frac{a_j}{q^{\alpha_j}}\right)\right)
\psi_0 \left(q^{\frac{2}{s}}\left(y_j-\frac{a'_j}{q^{\alpha_j}}\right)\right)
|x-y|^{-s}dxdy.$$

Then (\ref{acro}) becomes

\begin{equation}\label{submamad=2}
q^{-4}q^{\frac{8}{s}}
\sum \prod_{j=1}^2 
\psi_0 \left(\frac{a_j}{q^{\alpha_j}}\right) \psi_0 \left(\frac{a'_j}{q^{\alpha_j}}\right)I(a,a')
\end{equation}
where the sum is taken over $a\neq a'$, both in $\mathbb{Z}^2$.
 
We seek and upper bound for $|x-y|^{-s}$. 
Since $a\neq a'$, then either $a_1=a_1'$ and $a_2\neq a_2'$ (case 1), $a_1\neq a_1'$ and $a_2= a_2'$ (case 2), or $a_1\neq a_1'$ and $a_2\neq a_2'$ (case 3).
The first two cases are handled similarly and we omit the proof of case 2. 
\\

\textbf{Case 1:}
$$|x-y| \sim |x_1-y_1| + |x_2 -y_2| \geq q^{-\alpha_2} |a_2-a_2'|.$$
Now $$I(a,a') \lesssim  q^{-\frac{8}{s}}q^{s\alpha_2} |a_2-a_2'|^{-s}.$$ 

Thus,  (\ref{submamad=2})  is bounded above by
\begin{equation}\label{acrobatics}
q^{-4}q^{s\alpha_2}
\sum
\prod_{j=1}^2 \psi_0 \left(\frac{a_j}{q^{\alpha_j}}\right) \psi_0 \left(\frac{a'_j}{q^{\alpha_j}}\right)
|a_2-a_2'|^{-s} 
\end{equation}
where the sum is taken over $a_2 \neq a'_2$, both in $\mathbb{Z}$, and over $a_1= a'_1$ in $\mathbb{Z}$. 

Fix $a_2' \in \mathbb{Z}$ with $|a_2'|\le C q^{\alpha_2}$ and set
$$A_m = \{ a_2 \in \mathbb{Z}: 2^m \le |a_2-a_2'| < 2^{m+1} \}$$ where $0\le m \le \log{ \lceil Cq^{\alpha_2}\rceil}$. 
We use the fact that $\#(A_m) \sim 2^{m}$. 
We have
$$\sum_{a_2\in\mathbb{Z}}
\psi_0 \left(\frac{a_2}{q^{\alpha_2}}\right) 
|a_2-a_2'|^{-s} 
\sim \sum_m \sum_{a_2 \in A_m}   |a_2-a_2'|^{-s}   
\lesssim  \sum_m 2^{m(1-s)}
\sim 1.$$ 

Now, $$(\ref{acrobatics}) \le  q^{-4}q^{s\alpha_2} q^{\alpha_1} q^{\alpha_2}$$ which is bounded by 1 as $\alpha_2 \le \frac{2}{s}$ and $\alpha_1 +\alpha_2 =2$.

\vskip.125in 

\textbf{Case 3:} 
$$|x-y|^{-s} \geq q^{\frac{s(\alpha_1 + \alpha_2)}{2}}|a_1-a_2|^{-\frac{s}{2}}|a_1'-a_2'|^{-\frac{s}{2}}$$ 
to see that
(\ref{submamad=2})  is bounded above by
\begin{equation}\label{moreacrobatics}
q^{-4}q^{\frac{s(\alpha_1 + \alpha_2)}{2}}
\sum
\prod_{j=1}^2 \psi_0 \left(\frac{a_j}{q^{\alpha_j}}\right) \psi_0 \left(\frac{a'_j}{q^{\alpha_j}}\right)
|a_1-a_2|^{-\frac{s}{2}}|a_1'-a_2'|^{-\frac{s}{2}}
\end{equation}
where the sum is taken over $a_2 \neq a'_2$, both in $\mathbb{Z}$, and over $a_1= a'_1$ in $\mathbb{Z}$. For a fixed  $a_2\ \in \mathbb{Z}$ with $|a_2'|\le C q^{\alpha_2}$, let $A_m$ be as in case 1. 

Now 
$$ \sum_{a_2 \in \mathbb{Z}}  \psi_0 \left(\frac{a_2}{q^{\alpha_2}}\right)  |a_2-a_2'|^{-\frac{s}{2}}   
\sim \sum_m \sum_{a \in A_m}   |a_2-a_2'|^{-\frac{s}{2}}   
\lesssim  \sum_m 2^{m(1-\frac{s}{2})}
\sim q^{\alpha_2(1-\frac{s}{2})}.$$
Likewise, 
$$ \sum_{a_1 \in \mathbb{Z}}  \psi_0 \left(\frac{a_1}{q^{\alpha_1}}\right)  |a_1-a_1'|^{-\frac{s}{2}}   
\lesssim q^{\alpha_1(1-\frac{s}{2})}.$$

Therefore 
$$(\ref{moreacrobatics}) \lesssim q^{-4} q^{s\frac{\alpha_1 + \alpha_2}{2}}q^{\alpha_1 + \alpha_2}q^{\alpha_1(1-\frac{s}{2})} q^{\alpha_2(1-\frac{s}{2})} \lesssim 1.$$

\vskip.125in 

\textbf{The general (non-isotropic) case in all dimensions:} We are now ready to present the general case. We again break the sum over $a, a' \in \mathbb{Z}^d$ into the sum over the diagonal and the sum away from the diagonal. 
When  $a=a'$, (\ref{acro}) becomes

$$q^{-2d}q^{\frac{2d^2}{s}}
\sum_{a\in \mathbb{Z}^d}\prod_{j=1}^d\psi_0 \left(\frac{a_j}{q^{\alpha_j}}\right) 
\iint \psi_0 \left(q^{\frac{d}{s}}\left(x_j-\frac{a_j}{q^{\alpha_j}}\right)\right)\psi_0 \left(q^{\frac{d}{s}}\left(y_j-\frac{a_j}{q^{\alpha_j}}\right)\right)|x-y|^{-s}dxdy 
$$
and we use the same technique as above to show that this is bounded. 

In the case that $a\neq a'$, then $a_j\neq a'_j$ for at least one choice of $1\le j \le d$. Let $i\geq 1$ denote the number of coordinates for which $a_j\neq a'_j$. In the language of coding theory this means that the Hamming distance between $a$ and $a'$ is $i$. Choose a permutation of $\{1,..., d\}$, call it $\sigma$, such that $a_{\sigma(j)} \neq a'_{\sigma(j)}$ for $1\le j \le i$ and $a_{\sigma(j)} = a'_{\sigma(j)}$ for $i< j \le d$. 

Set $$I(a,a') = \iint \prod_{j=1}^d \psi_0 \left(q^{\frac{d}{s}}\left(x_{\sigma(j)}-\frac{a_{\sigma(j)}}{q^{\alpha_{\sigma(j)}}}\right)\right)
\psi_0 \left(q^{\frac{d}{s}}\left(y_{\sigma(j)}-\frac{a'_{\sigma(j)}}{q^{\alpha_{\sigma(j)}}}\right)\right)
|x-y|^{-s}dxdy.$$
Then (\ref{acro}) becomes

\begin{equation}\label{submama}
q^{-2d}q^{\frac{2d^2}{s}}
\sum \prod_{j=1}^d 
\psi_0 \left(\frac{a_{\sigma(j)}}{q^{\alpha_{\sigma(j)}}}\right) \psi_0 \left(\frac{a'_{\sigma(j)}}{q^{\alpha_{\sigma(j)}}}\right)I(a,a')
\end{equation}
where the sum is taken over $a_{\sigma(j)}\neq a'_{\sigma(j)}$ both in $\mathbb{Z}$, for $1\le j \le i$, and over $a_{\sigma(j)}= a'_{\sigma(j)}$  in $\mathbb{Z}$, for $i< j \le d$.

We observe that 
\begin{equation}\label{triangleineq} |x-y|^{-s} 
\le  i^{-s} \prod_{j=1}^i \left(\frac{|a_{\sigma(j)}- a'_{\sigma(j)}|}{q^{\alpha_{\sigma(j)}}}\right)^{-\frac{s}{i}}.\end{equation}
Certainly 
$$|x-y| \sim \sum_{j=1}^i|x_{\sigma(j)}- y_{\sigma(j)}| $$ 
and 
$$|x_{\sigma(j)}- y_{\sigma(j)}| \sim  \frac{|a_{\sigma(j)}- a'_{\sigma(j)}|}{q^{\alpha_{\sigma(j)}}}.$$

Combining these observations with the fact that the  arithmetic mean dominates the geometric mean \footnote{We simply mean the classical inequality ${(A_1 \cdot A_2 \cdot \dots \cdot A_n)}^{\frac{1}{n}} \leq \frac{A_1+A_2+\dots+A_n}{n}$.} verifies (\ref{triangleineq}).

Now $$I(a,a') \lesssim q^{-\frac{2d^2}{s}}\prod_{j=1}^i \left(\frac{|a_{\sigma(j)}- a'_{\sigma(j)}|}{q^{\alpha_{\sigma(j)}}}\right)^{-\frac{s}{i}}.$$

Thus,  (\ref{submama})  is bounded above by
\begin{equation}
q^{-2d}
\sum
\left( \prod_{j=1}^d \psi_0 \left(\frac{a_{\sigma(j)}}{q^{\alpha_{\sigma(j)}}}\right) \psi_0 \left(\frac{a'_{\sigma(j)}}{q^{\alpha_{\sigma(j)}}}\right)\right)
\left( \prod_{j=1}^i \left(\frac{|a_{\sigma(j)}- a'_{\sigma(j)}|}{q^{\alpha_{\sigma(j)}}}\right)^{-\frac{s}{i}}\right)
\end{equation}
where the sum is taken over $a_{\sigma(j)}\neq a'_{\sigma(j)}$ both in $\mathbb{Z}$, for $1\le j \le i$, and over $a_{\sigma(j)}= a'_{\sigma(j)}$ in $\mathbb{Z}$, for $i< j \le d$.\\

Summing in $a_{\sigma(j)}$, for $(i+1)\le j \le d$, we see that it suffices to show that
\begin{equation}\label{subsubmama}
q^{-2d}q^{\alpha_{\sigma(i+1)} + \cdots + \alpha_{\sigma(d)}}
\sum
 \prod_{j=1}^i \psi_0 \left(\frac{a_{\sigma(j)}}{q^{\alpha_{\sigma(j)}}}\right) \psi_0 \left(\frac{a'_{\sigma(j)}}{q^{\alpha_{\sigma(j)}}}\right)
 \left(\frac{|a_{\sigma(j)}- a'_{\sigma(j)}|}{q^{\alpha_{\sigma(j)}}}\right)^{-\frac{s}{i}}
\lesssim 1
\end{equation}
where the sum is taken over $a_{\sigma(j)}\neq a'_{\sigma(j)}$ both in $\mathbb{Z}$, for $1\le j \le i$.

For $1\le j \le i$, we show that
\begin{equation}\label{jsmall}
\sum_{\stackrel{a_{\sigma(j)}, a'_{\sigma(j)}\in \mathbb{Z}}{a_{\sigma(j)}\neq a'_{\sigma(j)}}}
\psi_0 \left(\frac{a_{\sigma(j)}}{q^{\alpha_{\sigma(j)}}}\right) \psi_0 \left(\frac{a'_{\sigma(j)}}{q^{\alpha_{\sigma(j)}}}\right)
\left(\frac{|a_{\sigma(j)}- a'_{\sigma(j)}|}{q^{\alpha_{\sigma(j)}}}\right)^{-\frac{s}{i}} \lesssim q^{\alpha_{\sigma(j)}(1+\frac{s}{i})} \text{, for $s\geq i$,} 
\end{equation} 
and that when $s<i$, the left hand side of (\ref{jsmall}) is bounded by $q^{2\alpha_{\sigma(j)}}$.

Taking (\ref{jsmall}) for granted, we may complete the proof of the lemma. Indeed, recalling that $\alpha_j \le \frac{d}{s}$, we conclude that for $s\geq i$ 
$$ (\ref{subsubmama})\lesssim
q^{-2d}q^{\alpha_{\sigma(i+1)} + \cdots + \alpha_{\sigma(d)}} \prod_{j=1}^i q^{\alpha_{\sigma(j)}(1+\frac{s}{i})} \le 1.$$

For  $s<i$, taking (\ref{jsmall}) for granted, we recall that $\sum_{j=1}^d\alpha_j=d$ to conclude that
$$(\ref{subsubmama}) \lesssim
q^{-2d}q^{\alpha_{\sigma(i+1)} + \cdots + \alpha_{\sigma(d)}}\prod_{j=1}^iq^{2\alpha_{\sigma(j)}} \le 1.$$

This completes the proof of the lemma modulo the proof of (\ref{jsmall}). \\

We now prove (\ref{jsmall}). 

Fix $a'_{\sigma(j)} \in \mathbb{Z}$ such that $a'_{\sigma(j)} \le C q^{\alpha_{\sigma(j)}} $. 
Set 
$$A_m = \{ a_{\sigma(j)}\in \mathbb{Z}: 2^m \le |a_{\sigma(j)}-a'_{\sigma(j)}| < 2^{m+1} \}$$ where $0\le m$ and $2^{m+1} = \lceil Cq^{\alpha_{\sigma(j)}}\rceil$. 

Then
$$\sum_{\stackrel{a_{\sigma(j)}\in \mathbb{Z}}{a_{\sigma(j)}\neq a'_{\sigma(j)}}}  \psi_0 \left(\frac{a_{\sigma(j)}}{q^{\alpha_{\sigma(j)}}}\right) |a_{\sigma(j)}- a'_{\sigma(j)}|^{-\frac{s}{i}} \lesssim \sum_m \sum_{A_m}  |a_{\sigma(j)}- a'_{\sigma(j)}|^{-\frac{s}{i}}
\lesssim \sum_m   2^{m(1-\frac{s}{i})}.$$
If $s\geq i$, then $\sum_m   2^{m(1-\frac{s}{i})} \lesssim 1$, and so
\begin{align*}
&\sum_{\stackrel{a_{\sigma(j)}, a'_{\sigma(j)}\in \mathbb{Z}}{a_{\sigma(j)}\neq a'_{\sigma(j)}}}
\psi_0 \left(\frac{a_{\sigma(j)}}{q^{\alpha_{\sigma(j)}}}\right) \psi_0 \left(\frac{a'_{\sigma(j)}}{q^{\alpha_{\sigma(j)}}}\right)
\left(\frac{|a_{\sigma(j)}- a'_{\sigma(j)}|}{q^{\alpha_{\sigma(j)}}}\right)^{-\frac{s}{i}}  
\lesssim q^{\alpha_{\sigma(j)}(1+\frac{s}{i})}.
\end{align*}

If $s<i$, then $\sum_m   2^{m(1-\frac{s}{i})} \lesssim q^{\alpha(1-\frac{s}{i})}$, and so
\begin{align*}
&\sum_{\stackrel{a_{\sigma(j)}, a'_{\sigma(j)}\in \mathbb{Z}}{a_{\sigma(j)}\neq a'_{\sigma(j)}}}
\psi_0 \left(\frac{a_{\sigma(j)}}{q^{\alpha_{\sigma(j)}}}\right) \psi_0 \left(\frac{a'_{\sigma(j)}}{q^{\alpha_{\sigma(j)}}}\right)
\left(\frac{|a_{\sigma(j)}- a'_{\sigma(j)}|}{q^{\alpha_{\sigma(j)}}}\right)^{-\frac{s}{i}}  
\lesssim q^{2\alpha_{\sigma(j)}}.
\end{align*}
This completes the proof of the lemma. 

\vskip.125in 

\subsection{Estimation of the second term in (\ref{heaven})} It remains to show that 
\begin{equation}\label{Tbound} \|(\widehat{T_q\mu_q}(\cdot)\times|\cdot|^{(d-s)/2}) \|_2 \lesssim 1,\end{equation}
when $s \ge \frac{d+1}{2}$. 
Since the Monge-Ampere determinant of $\phi$ does not vanish 
on the set 
$$\{(x,y): \phi(x,y)=t \}$$ for $t>0$, $\phi$ satisfies the Phong-Stein rotational curvature condition of Phong and Stein on this set (\cite{PS86}, \cite{PhSt91}), and thus
\begin{equation}\label{mapping} T_q: L^2(\mathbb{R}^d) \rightarrow L^2_{\frac{d-1}{2}}(\mathbb{R}^d)\end{equation}
with constants uniform in t and $q$. See also \cite{St93} and \cite{So93} for the background and a thorough description of these are related estimates. 

We shall deduce (\ref{Tbound}) from the following result. 
\begin{proposition}\label{energybonus} 
Let $f$ be a Schwartz class function with with finite $s-dimensional$ energy, as defined in Lemma \ref{energy}, for $\frac{d+1}{2}\le s<d$. Suppose that ${||f||}_1 \leq C$ for some uniform constant $C>0$. Then 
\begin{equation}\label{energybonusestimate}\|(\widehat{T_qf}(\cdot)\times|\cdot|^{(d-s)/2}) \|_2 \lesssim 1.\end{equation}
\end{proposition} 

\section{Proof of Proposition \ref{energybonus} }
We fix positive Schwartz class functions $\eta_0(\xi)$ supported in the ball $\{ |\xi| \leq 4 \}$ and $\eta(\xi)$ supported in the spherical shell 
$$\{ 1 < |\xi| < 4\} \ \text{with} \ \eta_j(\xi)=\eta(2^{-j}\xi), \ j \geq 1,$$and 
$$\eta_0(\xi) + \sum_{j=1}^{\infty} \eta_j(\xi) =1.$$

Define the Littlewood-Paley piece of $f$ (see e.g. (\cite{St93})), denoted by $f_j$ for $j \geq0$, by the relation 
$$ \widehat{f}_j(\xi)=\widehat{f}(\xi) \eta_j(\xi).$$ 

Now the left hand side of Proposition \ref{energybonus} can be written as

\begin{equation}\label{sum}
\sum_{j\geq0}\sum_{k\geq0} \int \widehat{T_qf_j}(\xi) \overline{\widehat{T_qf_k}}(\xi)  |\xi|^{d-s} d\xi.
\end{equation}

We handle the sums in (\ref{sum}) in three steps. First, we fix $j \geq0$ and consider the case when $j=k$. Second, we consider the more general scenario where $|j-k|\le 2L$. This second step follows as a simple consequence of the first. Finally, we handle the case when $|j-k|>2L$. Here $L$ is some positive number to be determined. 
\\

The proof of the following lemma is provided following the proof of Proposition \ref{energybonus} and can be found in  both \cite{EIT10} and \cite{IJL10}.
\begin{lemma}\label{estimatefar} Let $\eta_l$ be as above. Then for any $M>0$, there exists a constant $C_m>0$ and an $L= L_M >0$ so that 
$$|\widehat{T_qf_j}(\xi)| \eta_l(\xi) \le C_M 2^{-M(\max\{j,l\})} $$
whenever $|j-l|>L$. 
\end{lemma}


\vskip.125in

\textbf{Case 1 ($j=k$):}
We first establish that (\ref{sum}) holds when $j=k$. That is,
\begin{equation}\label{motivation}
\sum_{j\geq0} \int|\widehat{T_qf_j}(\xi)|^2|\xi|^{d-s}  d\xi\lesssim 1.\end{equation}

To see this, decompose the integral by writing the left hand side of (\ref{motivation}) as
\begin{equation}\label{jequalk}\sum_{j\geq0} \sum_{l\geq0} \int \eta_l(\xi) |\widehat{T_qf_j}(\xi)|^2|\xi|^{d-s}  d\xi.\end{equation}

Next, fix $j\geq 0$ and consider both the sum over $l\geq0$ such that $|j-l|\le L$ and the sum over $l>0$ such that $|j-l|>L$, where $L$ is some positive number to be determined. That is,

\begin{align*}
(\ref{jequalk})
=& \sum_{j\geq0}\sum_{\stackrel{l\geq0}{|j-l|\le L}} \int \eta_l(\xi) |\widehat{T_qf_j}(\xi)|^2|\xi|^{d-s}  d\xi
+ \sum_{j\geq 0}\sum_{\stackrel{l\geq0}{|j-l|> L}}\int \eta_l(\xi) |\widehat{T_qf_j}(\xi)|^2|\xi|^{d-s}  d\xi.
\\
=& \hspace{0.2in} I + \hspace{0.2in} II.
\\
\end{align*}

For $|j-l|\le L$, we use the support conditions for $\eta_l$ and the mapping properties of $T_q$ to write
\begin{align*}
I
&\sim \sum_{j\geq0} \sum_{\stackrel{l\geq0}{|j-l|\le L}} 2^{l(1-s)}\int \eta_l(\xi) |\widehat{T_qf_j}(\xi)|^2 |\xi|^{d-1}  d\xi
\\
&\lesssim \sum_{j\geq0} \sum_{\stackrel{l\geq0}{|j-l|\le L}} 2^{l(1-s)}\int  |\widehat{f_j}(\xi)|^2 d\xi
\\
&\le \sum_{j\geq0} \sum_{\stackrel{l\geq0}{|j-l|\le L}} 2^{(j-L)(1-s)}\int  |\widehat{f_j}(\xi)|^2 d\xi
\\
&\le (2L+1) 2^{-L(1-s)}\sum_{j\geq0}2^{j(1-s)}\int  |\widehat{f_j}(\xi)|^2 d\xi.
\\
\end{align*}

Since $(1-s)\le (s-d)$, when $\frac{d+1}{2}\le s <d$, and since $\widehat{f_j}(\xi)$ is supported where $|\xi| \sim 2^j$ then
\begin{align*}
I
&\le (2L+1) 2^{-L(1-s)}\sum_{j\geq0}2^{j(s-d)}\int  |\widehat{f_j}(\xi)|^2 d\xi.
\\
&\sim (2L+1) 2^{-L(1-s)}\sum_{j\geq0}\int  |\widehat{f_j}(\xi)|^2 |\xi|^{s-d} d\xi 
\\
&\sim (2L+1) 2^{-L(1-s)} \int  |\widehat{f}(\xi)|^2 |\xi|^{s-d} d\xi.
\\
\end{align*}
Finally, since $f$ has finite $s-dimensional$ energy when $\frac{d+1}{2}\le s <d$, then $ I \lesssim 1$.

To bound $II$, use Lemma \ref{estimatefar} to write
$$II \lesssim \sum_{j}\sum_{\stackrel{l\geq0}{|j-l|> L}}   2^{-M(\max\{j,l\})}   \int \eta_l(\xi) |\xi|^{d-s}  d\xi.$$
Since $\eta_l$ is compactly supported, we conclude that $II\lesssim 1$ thus finishing the proof of the first case. 
\vskip.125in

\textbf{Case 2 ($|j-k|\le 2L$):}
We now proceed to bounding (\ref{sum}) when $|j-k|\le 2L$.
We have
\begin{align*}
&\sum_{j\geq0}\sum_{\stackrel{k\geq 0}{|j-k|\le 2L}} \int \widehat{T_qf_j}(\xi) \overline{\widehat{T_qf_k}}(\xi)  |\xi|^{d-s} d\xi
\\
\le&\sum_{j\geq0}\sum_{\stackrel{k\geq 0}{|j-k|\le 2L}} \left(\int |\widehat{T_qf_j}(\xi) |^2 |\xi|^{d-s} d\xi\right)^{\frac{1}{2}}
\left(\int |\widehat{T_qf_k}(\xi)|^2 |\xi|^{d-s} d\xi\right)^{\frac{1}{2}}
\\
=& \sum_{|i|\le 2L} \sum_{j\geq0}\left(\int |\widehat{T_qf_j}(\xi) |^2 |\xi|^{d-s} d\xi\right)^{\frac{1}{2}}
\left(\int |\widehat{T_qf_{j+i}}(\xi)|^2 |\xi|^{d-s} d\xi\right)^{\frac{1}{2}}
\\
\le& \sum_{|i|\le 2L}  \left(\sum_{j\geq0}\int |\widehat{T_qf_j}(\xi) |^2 |\xi|^{d-s} d\xi\right)^{\frac{1}{2}}
\left(\sum_{j\geq0}\int |\widehat{T_qf_{j+i}}(\xi)|^2 |\xi|^{d-s} d\xi\right)^{\frac{1}{2}}
\\
\lesssim& 1,
\end{align*}
where we have applied the Cauchy Scwartz inequality twice and applied (\ref{motivation}).
\vskip.125in 

\textbf{Case 3 ($|j-k|>2L$):}
We now show that (\ref{sum}) is bounded when $|j-k| > 2L$.
We again decompose the integral by writing
$$\sum_{l\geq0} \sum_{j\geq0}\sum_{|j-k|> 2L} \int \eta_l(\xi) \widehat{T_qf_j}(\xi) \overline{\widehat{T_qf_k}}(\xi)  |\xi|^{d-s} d\xi.$$
Since $|j-k|>2L$, then either $|j-l|>L$ or $|k-l|>L$. 
Therefore, we may use Lemma \ref{estimatefar} to finish case 3. 

In more detail, assume $|j-l|>L$. Then 
$$\eta_l(\xi)| \widehat{T_qf_j}(\xi)| \lesssim 2^{-M\{j,l\}}.$$

If $|k-l|>L$, then 
$$\eta_l(\xi)| \widehat{T_qf_k}(\xi)| \lesssim 2^{-M\{j,l\}},$$ otherwise we use the observation that, for a fixed $l$, $\{k: |k-l|\le L\}$ is a finite set. 

This completes the proof of estimate (\ref{keyestimate}) and hence the proof of the theorem up to the proof of Lemma \ref{estimatefar}.

\section{Proof of Lemma \ref{sum}}
Recall from $(\ref{operator})$,
$$T_q f_j(x)= q^{\frac{d}{s}}\int_{\{y\in \mathbb{R}^d : |\phi(x,y)-1| \le q^{-\frac{d}{s}} \}}f_j(y)\psi(x,y) dy.$$

By Fourier inversion
$$T_q f_j(x)= \iiint_{\{y\in \mathbb{R}^d : |\phi(x,y)-1| \le q^{-\frac{d}{s}} \}} e^{i y \cdot \zeta} e^{is \cdot (\phi(x,y)-1)}  \psi(x,y)   \widehat{f_j}(\zeta) \widehat{\psi_0}(sq^{-\frac{d}{s}}) dy d\zeta ds,$$ 
where $\psi_0$ is a smooth compactly supported function centered at the origin and the change in the order of integration can be justified by Fubini's Theorem. 
\\

So
$$\widehat{T_q f_j}(\xi)= \iiiint_{\{y\in \mathbb{R}^d : |\phi(x,y)-1| \le q^{-\frac{d}{s}} \}} e^{-ix \cdot \xi} e^{iy \cdot \zeta}e^{is \cdot (\phi(x,y)-1)}  \psi(x,y)\widehat{f_j}(\zeta) \widehat{\psi_0}(sq^{-\frac{d}{s}}) dy d\zeta dsdx,$$
and therefore
\begin{align*} 
\widehat{T_q f_j}(\xi)\eta_l(\xi) 
=& \int   I_{jl}(\xi, \zeta, s)  \widehat{f}(\zeta) \widehat{\psi_0}(sq^{-\frac{d}{s}})d\zeta ds 
\\
\end{align*}
where 
\begin{equation}\label{oscillatory}I_{jl}(\xi, \zeta, s)  =\eta_l(\xi) \eta_j(\zeta) \iint_{\{y : |\phi(x,y)-1| \le q^{-\frac{d}{s}} \}} e^{i[s \cdot (\phi(x,y)-1) + y \cdot \zeta - x\cdot \xi]} \psi(x,y)dy dx .\end{equation}

Computing the critical points, $(x,y)$, of the phase function in (\ref{oscillatory}), we see that 
$$s \nabla_x\phi_l(x,y)=\xi  \text{ and } s \nabla_y \phi_l(x,y)=-\zeta.$$

The compact support of $\psi$ along with the non-zero gradient condition from (\ref{manifoldstructure}) implies that 
$$ \left| \nabla_x\phi_l(x,y)\right| \approx \left|\nabla_y\phi_l(x,y)\right| \approx 1.$$  
More precisely, the upper bound follows from smoothness and compact support. The lower bound follows from the fact that a continuous non-negative function achieves its minimum on a compact set. 
This minimum is not zero because of the condition (\ref{manifoldstructure}). 

It follows that 
\begin{equation} \label{orgasm} |\xi| \approx |\zeta| \end{equation} when we are near critical points. 
However, by comparing the support of $\eta_l$ with that of $\eta_j$ when $|j-l|>L$ we see that the integrand is supported away from critical points as  (\ref{orgasm}) no longer holds. 
This implies that for each non-critical point, $(x,y)$, either
\begin{equation}\label{good}
s \nabla_x\phi_l(x,y)\neq \xi \text{ or } s \nabla_y \phi_l(x,y)\neq-\zeta.
\end{equation}
Notice that this condition may vary with the choice of $(x,y)$. This will not, however, ultimately affect the argument due to the smoothness of $\phi$ and the presence of the compact function $\psi$ in the integrand. That is, we may restrict our attention to an open set containing a fixed non-critical point on which one of the equations holds, by restricting the support of $\psi$. Then we may repeat the argument over finitely many such open sets. 

Without loss of generality, assume that $l>j$.

Consider the case that $|s| >> |\xi|$ (i.e $|s| \geq c|\xi|$ with a sufficiently small constant $c>0$). We observe that, since $\left|\nabla_x\phi_l(x,y)\right| \approx 1 $, $\exists h$ so that $|\frac{\partial \phi}{\partial x_h}(x,y)|\approx 1$. 

It is immediate that $e^{-i x \cdot \xi} e^{is \cdot(\phi(x,y)-1)}$ is an eigenfunction of the differential operator 
$$L = \frac{1}{i(s \frac{\partial \phi}{\partial x_h} - \xi_h)} \frac{\partial}{\partial x_h}.$$

We will integrate by parts in (\ref{oscillatory}) using this operator.  The expression that we get after performing this procedure $M$ times is 
$$\left|I_{jl}(\xi, \zeta, s)\right| \lesssim \sup_{x,y} \left| s \frac{\partial \phi}{\partial x_h} - \xi_h \right|^{-M}.$$ 

\noindent Notice,
$$\left|s\frac{\partial \phi}{\partial x_h} - \xi_h \right| \gtrsim \left| \hskip.1cm \left|s \frac{\partial \phi}{\partial x_h}\right| - |\xi_h| \hskip.1cm \right| \approx |s|> |\xi|.$$  

\noindent So, $$|I_{jl}(\xi, \zeta, s)| \lesssim |\xi|^{-M} \lesssim 2^{-Ml} .$$
\\
In the case that $|s| << |\xi|$ (i.e $|s| \leq c|\xi|$ with a sufficiently small constant $c>0$), we observe that, since $|\xi| \sim 1 $, $\exists h'$ so that $|\xi_{h'}| \sim |\xi| \sim 1$. 
We notice that $e^{-i x \cdot \xi} e^{is \cdot(\phi(x,y)-1)}$ is an eigenfunction of the differential operator 
$$L = \frac{1}{i(s \frac{\partial \phi}{\partial x_{h'}} - \xi_{h'})} \frac{\partial}{\partial x_{h'}}$$

and we again integrate by parts in (\ref{oscillatory}) using this operator.  The expression that we get after performing this procedure $M$ times is 
$$\left|I_{jl}(\xi, \zeta, s)\right| \lesssim \sup_{x,y} \left| s \frac{\partial \phi}{\partial x_{h'}} - \xi_{h'} \right|^{-M}.$$ 

\noindent Notice,
$$\left|s\frac{\partial \phi}{\partial x_{h'}} - \xi_{h'} \right| \gtrsim \left| \hskip.1cm \left|s \frac{\partial \phi}{\partial x_{h'}}\right| - |\xi| \hskip.1cm \right| \approx  |\xi|$$  

\noindent and we again conclude that $$|I_{jl}(\xi, \zeta, s)| \lesssim |\xi|^{-M} \lesssim 2^{-Ml} .$$

Last, we consider the case when $|s| \sim |\xi| \sim 2^{l}$. By (\ref{good}), $\exists h$ such that 
$$ |s| \approx \left|s\nabla_x\phi_l(x,y)\right| \approx |s| \left| \frac{\partial \phi}{\partial y_{h}}\right|$$ and we notice that 
$e^{-i x \cdot \xi} e^{is \cdot(\phi(x,y)-1)}$ is an eigenfunction of the differential operator 
$$L = \frac{1}{i(s \frac{\partial \phi}{\partial y_{h}} - \zeta_{h})} \frac{\partial}{\partial y_{h}}.$$

The result once again follows by repeated integration by parts. This concludes the proof of the lemma.


\end{document}